\documentclass[letterpaper, 10pt]{article}
\usepackage{graphicx}
\usepackage{float}
\usepackage[inner=0.75in, outer=0.75in, top=0.75in, bottom=1in, paperheight=11in, paperwidth=8.5in]{geometry}

\usepackage{array}
\usepackage{rotating}
\usepackage{multirow}

\usepackage{mathtools}
\usepackage{shuffle}
\usepackage[table, svgnames, dvipsnames,xcdraw]{xcolor}
\usepackage[square, numbers, comma, sort&compress]{natbib}
\usepackage{verbatim}
\usepackage{bm}
\usepackage{amsmath}
\usepackage{amsfonts}
\usepackage{amssymb}
\usepackage{amsthm}

\usepackage{epstopdf}
\usepackage[utf8]{inputenc}
\usepackage[english]{babel}
\usepackage{color}
\usepackage{forest}

\usepackage{pdfpages}
\usepackage{pgfplots}
\usepackage{pgfplotstable}
\usepgfplotslibrary{patchplots}
\pgfplotsset{compat=1.15}

\usepackage{microtype}
\usepackage{url}
\usepackage{subfiles}
\usepackage{caption}
\usepackage{subcaption}
\usepackage{parskip}
\usepackage[section]{placeins}

\usepackage{tikz}
\usepackage{tikz-cd}
\usetikzlibrary{arrows,calc,shapes,decorations.pathreplacing}
\usetikzlibrary{decorations.markings}
\usetikzlibrary{trees,positioning,fit,shapes}
\usetikzlibrary{backgrounds}
\usetikzlibrary{patterns, positioning, arrows}

\usepackage{framed}

\colorlet{shadecolor}{blue!20}

\usepackage[obeyFinal,textwidth=0.55in]{todonotes}
\usepackage[hidelinks]{hyperref}
\usepackage[capitalize]{cleveref}

\usepackage{enumitem}

\theoremstyle{plain}
\newtheorem{thm}{Theorem}[section]
\theoremstyle{definition}

\begin{document}	
\graphicspath{{./figures/}}

\title{\bf On the application of the generating series for nonlinear systems with polynomial stiffness}	
\author{T.\ Gowdridge, N.\ Dervilis, K.\ Worden \\
        Dynamics Research Group, Department of Mechanical Engineering, University of Sheffield \\
        Mappin Street, Sheffield S1 3JD, UK
	   }
\date{}
\maketitle
\thispagestyle{empty}

\section*{Abstract}
Analytical solutions to nonlinear differential equations -- where they exist at all -- can often be very difficult to find. For example, Duffing’s equation for a 
system with cubic stiffness requires the use of elliptic functions in the exact solution. A system with general polynomial stiffness would be even more difficult to 
solve analytically, if such a solution was even to exist. Perturbation and series solutions are possible, but become increasingly demanding as the order of solution 
increases. This paper aims to revisit, present and discuss a geometric/algebraic method of determining system response which lends itself to automation. The method,
originally due to Fliess and co-workers, makes use of the generating series and shuffle product, mathematical ideas founded in differential geometry and abstract 
algebra. A family of nonlinear differential equations with polynomial stiffness is considered; the process of manipulating a series expansion into the generating 
series follows and is shown to provide a recursive schematic, which is amenable to computer algebra. The inverse Laplace-Borel transform is then applied to derive
a time-domain response. New solutions are presented for systems with general polynomial stiffness, both for deterministic and Gaussian white-noise excitation.

\textbf{Key words: Generating series; Nonlinear system modelling; Polynomial stiffness; Gaussian white noise.}

\section{Introduction}
\label{sec:intro}

One of the main problems in the discipline of nonlinear structural dynamics, is that the equations of motion concerned very very rarely admit closed-form exact 
solutions, and this usually forces a dependence on approximate solutions \cite{worden2019nonlinearity}. The most common approximation methods are series solutions,
often perturbation expansions in the coefficients of the nonlinear terms, which must be `small' in order that the series converge and that low-order truncations
are useful. An alternative series formulation is provided by the {\em Volterra series} \cite{volterra1930theory,barrett1963use}, which is not formally defined in
perturbative terms, but as a {\em functional series}. One of the attractive features of the Volterra series is that the generalised coefficients of the series have
physical interpretations; the coefficients are actually linear and nonlinear impulse responses, and their Fourier transforms can be regarded as Higher-dimensional 
Frequency Response Functions (HFRFs) \cite{worden2019nonlinearity}. Unfortunately, calculations with the Volterra series are very demanding in algebraic terms and 
rapidly become intractable (at least by hand), as the order of the expansion increases. One can bring computer algebra to bear on the problem, in order to automate 
calculations, but this does require a careful reformulation of the problem.

Fortunately, such a reformulation exists in the form of the {\em generating series}; this is an geometric-algebraic nonlinear system analysis developed in the 1980s,
by Michel Fliess and co-workers \cite{fliess1981fonctionnelles,fliess1982application,lamnabhi1982new,fliess1983algebraic,lamnabhi1986functional,lamnabhi1980application}.
The ground-breaking idea involved the representation of the Volterra series -- a sequence of high-dimensional integrals -- as a purely algebraic expansion. The cost
of the approach was that the expansion variables were {\em non-commutative}. The over-riding benefit of the approach was that the operations in the series expansion could
be implementated in computer algebra. Like the Laplace or Fourier transforms, the generating series offered a duality between time responses mediated by differential
equations and a purely algebraic approach based on polynomials in the transform variables. Once solutions were established in the transform domain, they could be 
taken back into the time domain using the {\em Laplace-Borel} transform \cite{lamnabhi1982nonlinearcomp,li2004formal}. Although the theory of the generating series is
extremely elegant, it sadly did not gain a great deal of traction in the structural dynamics community. However, with interest turning back towards HFRFs and concrete
calculations using the Volterra series \cite{worden2019nonlinearity}, it is arguably time to revisit the approach.

The modest objective of the current paper is to give an accelerated tutorial on the generating series and illustrate its use. A novel result presented here will be
from the analysis of a family of nonlinear systems with general polynomial terms; in the notation of \cite{fliess1981fonctionnelles}, the equation of motion is,

\begin{equation}
    \sum_{i=0}^n l_i\frac{d^i}{dt^i}y(t) + \sum^{m}_{i=2} \varepsilon_i y^i(t) = u(t)
\label{eq:diff system}
\end{equation}
where the first term is a general linear differential operator. For a Single-Degree-of-Freedom (SDOF) oscillator, $n = 2$ and the first and second derivatives
represent damping and inertia terms respectively. (The equation can always be scaled so that $l_n = 1$.) The second term on the LHS is the polynomial stiffness term.

With $n = 2$ and the order of nonlinearity $m = 3$, the equation becomes that of the asymmetric Duffing oscillator \cite{duffing1918erzwungene,worden2019nonlinearity},

\begin{equation}
    m \ddot{y} + c \dot{y} + k y + k_2 y^2 + k_3 y^3 = x(t)
\label{eq:aduff}
\end{equation}
in the standard notation where the overall (mass) scale has been restored; $y$ is the displacement (response) and $x$ is the force (excitation). Setting $k_2 = 0$ gives
the symmetric Duffing oscillator,

\begin{equation}
    m \ddot{y} + c \dot{y} + k y + k_3 y^3 = x(t)
\label{eq:sduff}
\end{equation}

Previous papers using the generating series showed results for systems with a single nonlinear term; in the current paper, results will be given for the asymmetric
quadratic-cubic equation (\ref{eq:aduff}). Responses for the system under both deterministic and random excitations will be considered. 

The layout of the paper is as follows: Section \ref{sec:background} will provide the basic terminology and definitions in order that one can motivate the generating
series. Section \ref{sec:consol} sketches the basic of related perturbation theory and shows how to construct a diagrammatic representation of functional expansions.
Section \ref{sec:gen_series} outlines the main ideas of the generating series in the context of the asymmetric Duffing oscillator system, and Section \ref{sec:response}
outlines how the series can be used to compute system responses. Section \ref{sec:noise} discusses how one calculates responses to white noise excitations and the 
paper then concludes.

\section{Background Theory and Definitions}
\label{sec:background}

\subsection{Algebraic Structure}

As mentioned above, the generating series is an expansion in non-commuting variables; this clearly means that the variables themselves are not standard real or complex
numbers and/or the product of the variables is not as standard. In fact both of these complexities are present, in the algebra of interest. The `basis' of the algebra
is provided by a set $X=\{x_0,x_1,...,x_n\}$ of symbols called the {\em alphabet}; the elements in $X$ are called {\em letters}. The alphabet $X$ generates a set $X^*$, 
which is called the {\em free monoid} over $X$, whose elements are sequences of the form $x_{j_v}...x_{j_0}$, and are called words
\cite{fliess1983algebraic, reutenauer1993freelie}. $X^*$ is thus the set of all words formed from the letters of the alphabet $X$. The `polymonial' terms of the 
generating series will be sets of words in some $X^*$, to be elaborated later. 

Having defined the variables in the algebra of interest, it remains to specify how they are multiplied together; this is by using the {\em shuffle product}. The shuffle 
product is a binary operator which represents the sum of all the interleaved products formed from a riffle shuffle over the two operands. The product is best explained
in terms of a number of basic identities, which are essential in working with the generating series:

\begin{enumerate}
    \item $1 \shuffle 1 = 1$
    \item $1 \shuffle w = w \shuffle 1 = w$
    \item $x_jw \shuffle x_j'w' = x_j(w \shuffle x_j'w') + x_j'(x_jw \shuffle w')$
    \item $x^k \shuffle x^{n-k} = \displaystyle {n \choose k} x^n$
\end{enumerate}

Of these identities, the most useful is arguably the third which provides a recursive means of evaluating the product when the words concerned are specified in terms of
their constituent letters. An example of the product showing the basic riffle nature is, 

\begin{displaymath}
            ab \shuffle cd = a(b \shuffle cd)+ c(ab \shuffle d) = abcd + acbd + acdb + cabd + cadb + cdab
\end{displaymath}

In the shuffle product, the terms inside each argument do not change their order. For example, in the third case listed above, $x_j$ always appears before $w$, and $x_j'$ 
always before $w'$, where $\{x_j,x_j',...\} \in X$ and $\{w,w',...\} \in X^*$ \cite{reutenauer1993freelie}. This recursive definition for the shuffle product of two 
generating series is discussed and explored in greater detail in section \ref{sec:shuff gen}. The recursion for the shuffle product is completed when at least one 
argument is reduced to the identity element (above in cases 1 and 2).

While this algebra may seem rather strange and counter-intuitive, it will be shown that it arise naturally in automating the transition between the nonlinear differential
equation of interest and the corresponding Volterra series and thence to the generating series. The algebra arises because of the presence of {\em iterated integrals}. 

\subsection{Iterated Integrals}
\label{subsec:iterated}

There is nothing mysterious about iterated integrals, they are simply multiple integrals with the individual integrations carried out in a prescribed order. It is 
well-known that changing the order of integration in a multiple integral will change the integrand, and this is what will cause non-commutativity here. As a matter of
notation, iterated integrals will be denoted here like,

\begin{equation}
    \int_0^t d{x_n}...d{x_0}
\label{eq:it_int}
\end{equation}
and the convention will be that individual integrations will work inwards from the right in terms of variables; here the first integral will be over $x_0$ and the last
will be over $x_n$. In this example, all the limits on the integrals are the same and so the integrals are represented by a single symbol; in the general case, each
integral would have its own symbol and limits and these would also be traversed working inwards, this time from the left. In the generating series algebra, single 
integrals will correspond to the letters and multiple integrals to the words of the corresponding free monoid.

It is possible to see how the shuffle product might arise on such an algebra; consider a product of two iterated integrals,

\begin{equation}
    \Big( \int_0^t d\xi_{j_{v}}...d\xi_{j_{0}} \Big) \Big( \int_0^t d\xi_{k_{\mu}}...d\xi_{k_{0}} \Big)
\label{eq:it_prod}
\end{equation}

After a certain amount of standard calculus, one finds that integration by parts results in,

\begin{multline}
\int_0^t d\xi_{j_v}(\tau)\Big[\Big(\int_0^\tau d\xi_{j_{v-1}}...d\xi_{j_{0}}\Big) \Big(\int_0^\tau d\xi_{k_{\mu}}...d\xi_{k_{0}}\Big)\Big]
+\int_0^t d\xi_{k_\mu}(\tau)\Big[\Big(\int_0^\tau d\xi_{j_{v}}...d\xi_{j_{0}}\Big) \Big(\int_0^\tau d\xi_{k_{\mu-1}}...d\xi_{k_{0}}\Big)\Big]
\end{multline}
and this is very suggestive of the relation for the shuffle product,

\begin{equation}
    (xw)\shuffle(x'w')=x[w\shuffle(x'w')]+x'[(xw)\shuffle w']
\label{eq:shuffle def}
\end{equation}

In fact, this is evidence of the correspondence with nonlinear differential equations; the nonlinear terms in such equations engender products of iterated
integrals which map to shuffle products in the algebra of the generating series. In the language of \cite{fliess1982application, ree1958lie}:

\begin{thm}
\label{theorem}
The product of two analytic causal functionals is again an analytic causal functional of the same kind, the generating power series of which is the shuffle product of the two generating power series.
Formally, this represents
\begin{equation}
    y_1 \times y_2 \Leftrightarrow g_1 \shuffle g_2
\end{equation}
\end{thm}

This theorem can be extended to higher-order products of terms. In the differential equation, terms of the form $y^n$, interpreted as $y\times ... \times y$, $n$ times,
map directly to `powers' in the generating series where the product is the shuffle. The simplest way to demonstrate this will be via the concrete examples to be pursued
shortly.

\subsection{Volterra Series}
\label{subsec:volterra}

As mentioned in the introduction, another key ingredient in methodology here is the {\em Volterra series}; this is essentially a functional Taylor series for the 
response of a nonlinear input-output system \cite{volterra1930theory,barrett1963use}. The series generalises the Duhamel integral for a linear system 
$x(t) \longrightarrow y(t)$, given by \cite{worden2019nonlinearity},

\begin{equation}
    y(t) = \int^{+\infty}_{-\infty} h(\tau)x(t-\tau) d\tau
\label{eq:duham}
\end{equation}
where $h(\tau)$ represents the impulse response of the system. For a nonlinear system, one obtains instead an infinite series,

\begin{equation}
    y(t) = y_0(t) + y_1(t) + y_2(t) + \ldots + y_i(t) + \ldots
\label{eq:volt}
\end{equation}
where the general term is,

\begin{equation}
     y_i(t) = \int_{-\infty}^{+\infty} \ldots \int_{-\infty}^{+\infty} h_i(\tau_1,...,\tau_i) x(t - \tau_1) \ldots x(t - \tau_i) d\tau_1 \ldots d\tau_i
\label{eq:kern}
\end{equation}
which is of course, an iterated integral. The functions $h_i(\tau_1,\ldots,\tau_i)$ are the `coefficients' in the functional expansion and have a direct physical 
interpretation as higher-dimensional impulse response functions \cite{worden2019nonlinearity}; they are termed {\em Volterra kernels}. Clearly, the problem of 
establishing a Volterra series is that of determining the kernels for a given nonlinear system. The multi-dimensional Fourier transforms also have a physical 
interpretation as higher-dimensional frequency response functions (HFRFs). One way to find the Volterra kernels is by determining the HFRFs directly and then using
an inverse Fourier transform; although this sounds rather indirect, it is possible because the HFRFs can be found from the nonlinear equations of motion by 
{\em harmonic probing} as introduced, in \cite{bedrosian1971output}.

In fact, a variant of the strategy just described will be used in this paper to find the Volterra kernels and system responses. Rather than computing objects in the 
Fourier domain and using an (inverse) Fourier transform to find time-domain objects, the idea will be to compute objects in the algebra of the generating series and 
transform back; the relevant transform is called the inverse {\em Laplace-Borel transform}. Like the Laplace transform, the forward map is easier to find than the 
inverse, so the usual means of inversion is to use a table of inverse Laplace-Borel transforms \cite{lamnabhi1982new,lamnabhi1986functional,lubbock1969multidimensional}.

\section{The Consolidated Expansion and a Diagrammatic Representation}
\label{sec:consol}

As mentioned earlier, the best way of illustrating the difficult concepts here is via concrete examples. With this in mind, this section will single out the 
asymmetric Duffing oscillator (equation (\ref{eq:aduff})) as the system of interest. The system is actually of considerable practical interest as it represents the
lowest-order approximation to a general SDOF system with both odd and even nonlinearities. Before proceeding, it is important to note that the form of equation can be
simplified without losing generality. By scaling the independent ($t$) and dependent variables ($x,y$), the number of parameters can be reduced, so that the equation 
becomes,

\begin{equation}
    \ddot{y} + a\dot{y} + y + \varepsilon_1 y^2 + \varepsilon_2 y^3 = x(t)
\label{eq:saduff}
\end{equation}
and it is this form which is considered from now on, in order to simplify the notation and algebra. Note that there are parameters associated with each of the nonlinear
terms, $\varepsilon_1$ (quadratic) and $\varepsilon_2$ (cubic); in a standard perturbation approach, these would be the expansion parameters, and this will also be the 
case here. It is useful at this point to look at the perturbation approach, even if it will not be pursued directly here; the formulation will show the complexity of 
the problem and also allow the construction of a useful and intuitive diagrammatic representation.

\subsection{Developing the Diagrammatic Representation}
\label{subsec:diag rep}

The first stage in the analysis here is to pass to the frequency domain via the Fourier transform. The standard operations on equation (\ref{eq:saduff}) yield, via the 
convolution theorem,

\begin{equation}
    (1 + i\alpha \omega - \omega^2) Y(\omega) + \varepsilon_1 (Y*Y)(\omega) + \varepsilon_2 (Y*Y*Y)(\omega)  
\label{eq:fsaduff}
\end{equation}
where,

\begin{equation}
    (Y*Y)(\omega) = {\int_{-\infty}^{\infty}} Y(\omega - \Omega) Y(\Omega) d \Omega
\label{eq:conv2}
\end{equation}
and,

\begin{equation}
    (Y*Y*Y)(\omega) = {\int_{-\infty}^{\infty}} {\int_{-\infty}^{\infty}} Y(\omega - \Omega) Y(\Omega - \Omega') Y(\Omega') d\Omega d\Omega'
\label{eq:conv3}
\end{equation}

These latter expressions are cumbersome and unsymmetrical and can be rewritten as symmetrical integrals,

\begin{equation}
    (Y*Y)(\omega) = {\int_{-\infty}^{\infty}} {\int_{-\infty}^{\infty}} Y(\omega_1) Y(\omega_2) \delta(\omega - \omega_1 - \omega_2) d \omega_1 d\omega_2
\label{eq:sconv2}
\end{equation}
and,

\begin{equation}
    (Y*Y*Y)(\omega) = {\int_{-\infty}^{\infty}} {\int_{-\infty}^{\infty}} {\int_{-\infty}^{\infty}} Y(\omega_1) Y(\omega_2) Y(\omega_3) 
    \delta(\omega - \omega_1 - \omega_2 - \omega_3) d\omega_1 d\omega_2 d\omega_3
\label{eq:sconv3}
\end{equation}

With these modifications, and a little more rearrangement, equation (\ref{eq:fsaduff}) becomes,

\begin{equation}
    Y(\omega) = H(\omega) X(\omega) - \varepsilon_1 H(\omega) {\int_{-\infty}^{\infty}} Y(\omega_1) Y(\omega_2) d\mu_2 - \varepsilon_2 
    H(\omega) {\int_{-\infty}^{\infty}} Y(\omega_1) Y(\omega_2) Y(\omega_3) d\mu_3
\label{eq:sfsaduff}
\end{equation}
where the integral signs have been coalesced and the measures are $d\mu_2 = \delta(\omega - \omega_1 - \omega_2) d \omega_1 d\omega_2$ and 
$d\mu_3 = \delta(\omega - \omega_1 - \omega_2 - \omega_3) d \omega_1 d\omega_2 d\omega_3$. Furthermore $H(\omega) = 1/(1 + i\alpha \omega - \omega^2)$, which is the 
standard FRF of the underlying linear system ($\varepsilon_1 = \varepsilon_2 = 0$). Note that the equation is {\em recursive}: i.e.\ $Y(\omega)$ is expressed as a 
nonlinear function of itself. In the case that $\varepsilon_1$ and $\varepsilon_2$ were small perturbation parameters, the equation could be used to compute 
$Y(\omega)$ iteratively, starting from the response of the linear system. With this observation in mind, it makes sense to compare the result with the actual
two-parameter perturbation expansion. In the time domain one has,

\begin{equation}
    y(t) = y_{00}(t) + \varepsilon_1 y_{10}(t) + \varepsilon_2 y_{01}(t) + \varepsilon_1^2 y_{20}(t) + \varepsilon_1^1 \varepsilon_2^1 y_{11}(t) + \ldots = \sum^\infty_{j=0}\sum^j_{i=0}\varepsilon^i_1\varepsilon^{j-i}_{2}y_{i,j-i}(t)
\label{eq:2nl_tpert}
\end{equation}
and the Fourier transform is,

\begin{equation}
     Y(\omega)=Y_{00}(\omega)+\varepsilon_1 Y_{10}(\omega) + \varepsilon_2 Y_{01}(\omega) + \varepsilon_1^2 Y_{20}(\omega) + \varepsilon_1^1 \varepsilon_2^1 Y_{11}(\omega) + ... = \sum^\infty_{j=0}\sum^j_{i=0}\varepsilon^i_1\varepsilon^{j-i}_{2}Y_{i,j-i}(\omega)
\label{eq:2nl fpert}
\end{equation}
with the obvious notation.

Now equating equations (\ref{eq:2nl fpert}) and (\ref{eq:sfsaduff}) at each level of perturbation, one obtains an infinite sequence of equations; the first nine,
corresponding to a truncation at third nonlinear order, are:

\begin{equation}
\begin{array}{lllll}
    O(\varepsilon_1^0 \varepsilon_2^0) & :~~~~~ &  Y_{00}(\omega) & = & H(\omega)X(\omega) \nonumber \\
    O(\varepsilon_1^1 \varepsilon_2^0) & :~~~~~ &  Y_{10}(\omega) & = & - \varepsilon_1 H(\omega) \int Y_{00}(\omega_1) Y_{00}(\omega_2) d\mu_2 \nonumber \\
    O(\varepsilon_1^0 \varepsilon_2^1) & :~~~~~ &  Y_{01}(\omega) & = & - \varepsilon_2 H(\omega)
    \int Y_{00}(\omega_1) Y_{00}(\omega_2) Y_{00}(\omega_3) d\mu_3 \nonumber \\
    O(\varepsilon_1^2 \varepsilon_2^0) & :~~~~~ & Y_{20}(\omega) & = & - \varepsilon_1 H(\omega) \int 2 Y_{00}(\omega_1) Y_{10}(\omega_2) d\mu_2 \nonumber \\
    O(\varepsilon_1^1 \varepsilon_2^1) & :~~~~~ & Y_{11}(\omega) & = & - \varepsilon_1 H(\omega) \int 2 Y_{00}(\omega_1) Y_{01}(\omega_2) d\mu_2 
    - \varepsilon_2 H(\omega) \int 3Y_{00}(\omega_1) Y_{00}(\omega_2) Y_{10}(\omega_3) d\mu_3 \nonumber \\
    O(\varepsilon_1^0 \varepsilon_2^2) & :~~~~~ & Y_{02}(\omega) & = & - \varepsilon_2 H(\omega) \int 3Y_{00}(\omega_1) Y_{00}(\omega_2) Y_{01}(\omega_3) d\mu_3 \nonumber \\
    O(\varepsilon_1^3 \varepsilon_2^0) & :~~~~~ & Y_{30}(\omega) & = & - \varepsilon_1 H(\omega) 
    \int [2Y_{00}(\omega_1) Y_{20}(\omega_2) + Y_{10}(\omega_1) Y_{10}(\omega_2)] d\mu_2 \nonumber \\
    O(\varepsilon_1^2 \varepsilon_2^1) & :~~~~~ & Y_{21}(\omega) & = & - \varepsilon_1 H(\omega) 
    \int [2Y_{00}(\omega_1) Y_{11}(\omega_2) + 2Y_{10}(\omega_1) Y_{01}(\omega_2)] d\mu_2 \nonumber \\ 
    ~ & ~ & ~ & ~ & - \varepsilon_2 H(\omega) \int [3Y_{00}(\omega_1) Y_{00}(\omega_2) Y_{20}(\omega_3) + 3Y_{00}(\omega_1) Y_{10}(\omega_2) Y_{10}(\omega_3)] 
    d\mu_3 \nonumber \\
    O(\varepsilon_1^1 \varepsilon_2^2) & :~~~~~ & Y_{12}(\omega) & = & - \varepsilon_1 H(\omega) 
    \int [2Y_{00}(\omega_1) Y_{02}(\omega_2) + Y_{01}(\omega_1) Y_{01}(\omega_2)] d\mu_2 \nonumber \\
    ~ & ~ & ~ & ~ & - \varepsilon_2 H(\omega) \int [6Y_{00}(\omega_1) Y_{10}(\omega_2) Y_{01}(\omega_3) + 3 Y_{00}(\omega_1) Y_{00}(\omega_2) Y_{11}(\omega_3)]
    d\mu_3 \nonumber \\ 
    O(\varepsilon_1^0 \varepsilon_2^3) & :~~~~~ &   Y_{03}(\omega) & = & - \varepsilon_2 H(\omega) \int [3Y_{00}(\omega_1) Y_{00}(\omega_2) Y_{02}(\omega_3) + 3Y_{00}(\omega_1) Y_{01}(\omega_2) Y_{01}(\omega_3)] d\mu_3 
\end{array}
\end{equation}
where the range of each integral is $(-\infty,\infty)$.

The expansion has been carried so far in order to show the contribution of multiple cross terms $\varepsilon_1^i \varepsilon_2^j$. Setting $\varepsilon_1 = 0$ (resp.\
$\varepsilon_2 = 0$) or collecting only the terms corresponding to $\varepsilon_1^i \varepsilon_2^0$ (resp.\ $\varepsilon_1^0 \varepsilon_2^j$) generates the expansion
for a lone quadratic (resp.\ cubic) nonlinearity. As in all perturbation expansions, each term is computable from previously evaluated terms, although the effort quickly
becomes large. In fact, one can see that $Y_{i,j}=\varepsilon_1 Y_{i-1,j} + \varepsilon_2 Y_{i,j-1}$. Although the algebraic representation provided here -- referred to 
as the {\em consolidated expansion} in \cite{morton1970} -- is cumbersome, the authors of that reference proposed a diagrammatic representation analogous to the Feynman
representation of perturbation expansions in quantum field theory \cite{feynman1948space}. The representation was also adopted in \cite{lamnabhi1986functional}; however,
it appears to have only been applied in the case of a single nonlinearity in previous work.

\begin{figure}[ht]
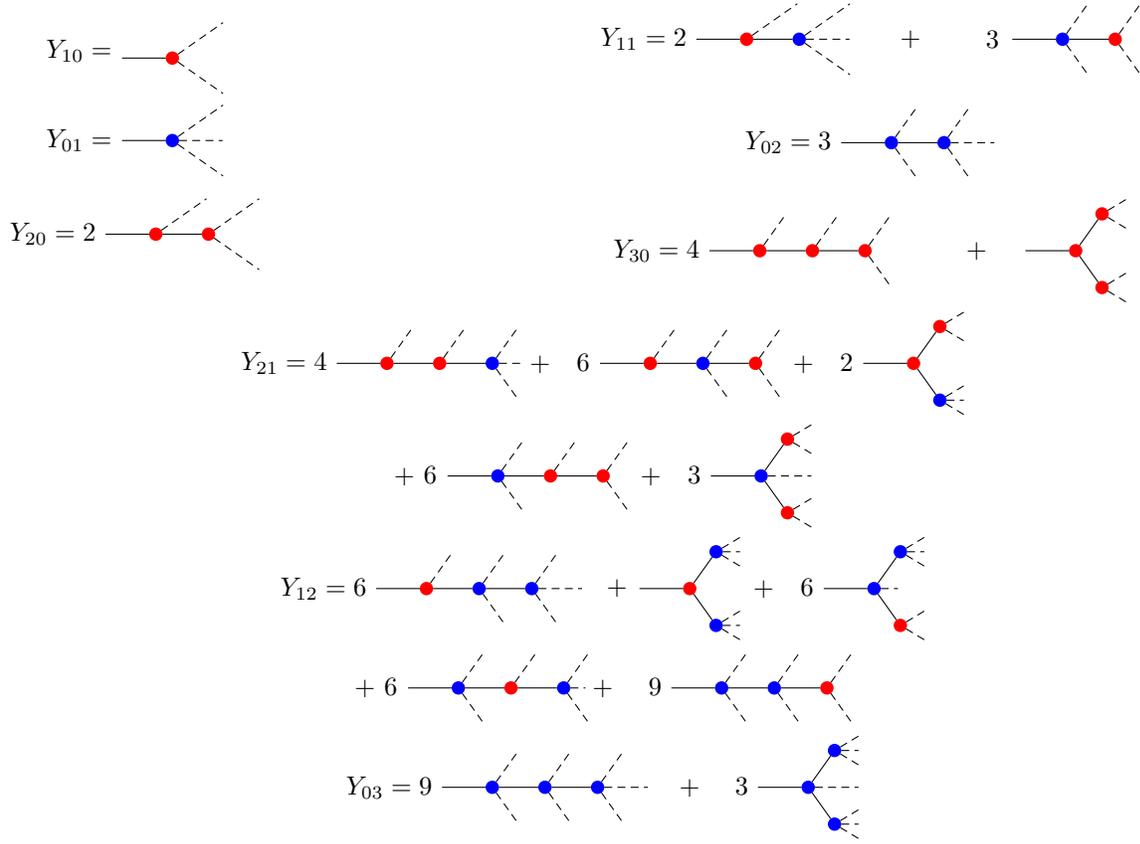

\centering
\include{figures/Figure_1}
\caption{Diagrammatic representation of first nine terms of the consolidated expansion for the quadratic-cubic oscillator.}
\label{fig:cons_exp}
\end{figure}

The correspondence between the algebraic expansion and the diagrammatic form is encoded in a set of rules, which allow each term to be represented by a tree diagram. In
the case of a quadratic-cubic system, the expansion is depicted in Figure \ref{fig:cons_exp}; the conventions are:

\begin{tikzpicture}[every node/.style={draw=black, fill=black, thick, circle, inner sep=1.5pt}]
\node[draw=white, inner sep=0pt] at (0,0) (no) {};
\node[draw=white, inner sep=0pt] at (0,0.1) (a) {};
\node[draw=white, inner sep=0pt] at (1,0.1) (b) {};
\draw (a) -- (b);
\end{tikzpicture}
corresponds to multiplication by $H(\omega)$.
          
\begin{tikzpicture}[every node/.style={draw=black, fill=black, thick, circle, inner sep=1.5pt}]
\node[draw=white, inner sep=0pt] at (0,0) (no) {};
\node[draw=white, inner sep=0pt] at (0,0.1) (a) {};
\node[draw=white, inner sep=0pt] at (1,0.1) (b) {}; 
\draw (a) [densely dashed]-- (b);
\end{tikzpicture} 
corresponds to multiplication by $Y_0(\omega_i)$, where the subscript $i$ can equal any positive integer.

\begin{tikzpicture}[every node/.style={draw=black, fill=black, thick, circle, inner sep=1.5pt}]
\node[draw=white, inner sep=0pt] at (0,0) (no) {};
\node[draw=white, inner sep=0pt] at (1,0.1) (b) {};
\node[draw=red, fill=red] at (0.5, 0.1) (l) {}; 
\end{tikzpicture} 
corresponds to multiplication by $\varepsilon_1$.

\begin{tikzpicture}[every node/.style={draw=black, fill=black, thick, circle, inner sep=1.5pt}]
\node[draw=white, inner sep=0pt] at (0,0) (no) {};
\node[draw=white, inner sep=0pt] at (1,0.1) (b) {};
\node[draw=blue, fill=blue] at (0.5, 0.1) (l) {}; 
\end{tikzpicture} 
corresponds to multiplication by $\varepsilon_2$

The rules associated with individual terms $Y_{ij}$ are:

\begin{description}
    \item[Rule 1:] The tree(s) will have $i$ vertices with 3 incident branches and $j$ vertices with 4 incident branches.
    \item[Rule 2:] There will be $i+j$ nodes and $i+j$ solid lines in the tree.
    \item[Rule 3:] A tree will have $i+2j+1$ dashed branches.
    \item[Rule 4:] Any two distinct vertices are connected by a single path.
    \item[Rule 5:] Frequency is conserved at a vertex. The sums of the frequencies either side of a vertex are equal.
\end{description}

The diagrammatic representation does not add anything to the algebraic expansion; however, exactly as it does in quantum field theory, it helps considerably in doing 
calculations by hand. Up to this point, the analysis has not really strayed beyond classical theory -- although Fliess did draw upon the diagrammatic representation in 
\cite{lamnabhi1986functional}. The generating series proper will be introduced in the next section via its calculation for the asymmetric Duffing oscillator.

\section{Generating Series for the Asymmetric Duffing Oscillator} 
\label{sec:gen_series}

As discussed above, the analysis will concentrate on the system specified in equation (\ref{eq:diff system}) at first, but then specialise to the asymmetric 
(quadratic-cubic) Duffing oscillator. The first stage in the analysis is to manipulate the equation into integral form; one integrates $n$ times in order to remove all
the differential operators, with result,

\begin{multline}
     y(t) + l_{n-1} \int^t_0y(\tau_1)+...+l_0\int^t_0d\tau_n\int^{\tau_n}_0d\tau_{n-1}...\int^{\tau_2}_0y(\tau_1)d\tau_1 + \\
     \sum^m_{i=2} \varepsilon_i \int_0^td\tau_n...\int^{\tau_2}_0y(\tau_1)d\tau_1=\int_0^td\tau_n...\int_0^{\tau_2}u(\tau_1)d\tau_1
\label{eq:gen_int}
\end{multline}
if one assumes that all relevant initial conditions are zero. The equation of motion now consists of iterated integrals and can be converted into the generating series 
domain. The symbol $g$ will denote the generating series associated with $y(t)$ here; the two basis letters associated with the free monoid will be denoted: 
$x_0=\int_0^t y(\tau)d\tau$ and $x_1=\int_0^t x(\tau)$. The rules for the transformation follow from Fliess' {\em fundamental formula} and the Peano-Baker formula
\cite{gantmakher2000theory}, as detailed in \cite{fliess1981fonctionnelles}. The main formal rules for the transformation are as follows:

\begin{enumerate}
\item The transform acts on linear combinations linearly.
\item Just as differentiation in the Fourier transform is represented by pre-multiplication by $i \omega$, integration in the $g$-domain is represented by 
      pre-multiplication by the word $x_0$.
\item $n^{th}$ powers of $y$ will transform to $n$-fold shuffle produces of $g$. 
\end{enumerate}

Applying these rules to equation (\ref{eq:gen_int}) yields,

\begin{equation}
    g + l_{n-1} x_0 g + \ldots + l_1 x_0^{n-1} g + 
    x_0^n (\varepsilon_2 g \shuffle g + \ldots + \varepsilon_m \underbrace{g\shuffle \ldots \shuffle g}_{m \text{ times}}) = x_0^{n-1} x_1
\label{eq:gen_poly}
\end{equation}
or the more compact form,

\begin{equation}
    (1 + \sum^{n-1}_{j=0} l_j x_0^{n-j}) g + x_0^n \sum^m_{i=2} \varepsilon_i g^{\shuffle i} = x_0^{n-1}x_1
\label{eq:cgen_poly}
\end{equation}
where $g^{\shuffle i} = g \shuffle \ldots \shuffle g$, $i$ times. 

Equation (\ref{eq:cgen_poly}) can be further simplified by factorising the the polynomial in $x_0$ that multiplies $g$ as follows,

\begin{equation}
    1 + \sum^{n-1}_{j=0} l_j x_0^{n-j} = \prod^p_{i=0} (1 - a_i x_0)^{\alpha_i}, ~~~ \alpha_1+\alpha_2+...+\alpha_p=n
\end{equation}

One can now formally write equation (\ref{eq:cgen_poly}) as,

\begin{equation}
    g = g_0 + \frac{x_0^n \sum^m_{i=2} \varepsilon_i g^{\shuffle i}}{\prod^p_{i=0} (1 - a_i x_0)^{\alpha_i}}
    \label{eq:gen_form}
\end{equation}
where,

\begin{equation}
   g_0 = \frac{x_0^{n-1} x_1}{\prod^p_{i=0} (1 - a_i x_0)^{\alpha_i}}
   \label{eq:def_g0}
\end{equation}
is seen to be the generating series representation of the underlying linear system. 

The solution to equation (\ref{eq:gen_poly}) can now be constructed recursively; starting with $g_0$, one computes,

\begin{equation}
    g_{i+1} = - \frac{x_0^{n} \times \sum^m_{j=2} \varepsilon_i \sum_{\nu_1+...+\nu_j = i} g_{\nu_1} \shuffle \ldots \shuffle g_{\nu_j}}
    {\prod^p_{i=0} (1 - a_i x_0)^{\alpha_i}}
\label{eq:gen_rec}
\end{equation}
and the representation of the full nonlinear system response is then,

\begin{equation}
    g = g_0 + g_1 + \ldots + g_i
\label{eq:gen_sum}
\end{equation}

At this point, it is important to recall that the algebra of the generating series is not commutative, so objects like equation (\ref{eq:def_g0}) are actually 
ambiguous. A careful analysis reveals that $g_0$, actually takes the form $R_1(x_0)x_{i_1}R_2(x_0)x_{i_2} \ldots x_{i_p}R_p(x_0)$ where $R_j(x_0)$ represents a rational 
fraction and $\{i_1, i_2, ... , i_p\} \in \{0,1\}$ \cite{fliess1983algebraic}. The quotient in the recursive scheme of equation (\ref{eq:gen_rec} is of a similar form; 
meaning that all the successive $g_i$ terms will also be of this form. In fact, the general form of the terms of interest can be written,

\begin{equation}
    \frac{1}{(1 - a_0 x_0)^{\alpha_0}} x_1 \frac{1}{(1 - a_1 x_0)^{\alpha_1}} x_1 \ldots x_1 \frac{1}{(1 - a_p x_0)^{\alpha_p}}
    \label{eq:gen_fact}
\end{equation}

Expressions of this type can be simplifying by using the identity,

\begin{equation}
    \frac{1}{(1 - a x_0)^{\alpha}} = \frac{1}{(1 - a x_0)^{\alpha-1}} + \frac{ax_0}{(1 - a x_0)(1 - a x_0)^{\alpha-1}}
    \label{eq:red_exp}
\end{equation}
and decomposing as partial fractions. In this way, by repeated application of the identity, all the exponents in the denominators can ultimately be reduced to unity, and 
the general object of interest becomes,

\begin{equation}
    \frac{1}{(1 - a_0 x_0)} x_{i_1} \frac{1}{(1 - a_1 x_0)} x_{i_2} \ldots x_{i_p} \frac{1}{(1 - a_p x_0)}
    \label{eq:gen_fact2}
\end{equation}

The general problem of computing shuffle products is thus reduced to that of computing shuffle produces of terms like that above. Summation of such terms is not an issue 
as the shuffle product is distributive over addition \cite{reutenauer1993freelie}.

\subsection{Shuffle Product of Series Terms} 
\label{sec:shuff gen}

Shuffle products of terms of the specific form in (\ref{eq:gen_fact2}) can now be considered in more detail. Suppose the two terms of interest are,

\begin{equation}
    g^p_1 = \frac{1}{1 - b_0 x_0} x_{i_1} \frac{1}{1 - b_1 x_0} x_{i_2} \ldots x_{i_p} \frac{1}{1 - b_p x_0} = g^{p-1}_1 \frac{x_{i_p}}{1 - b_p x_0}
\label{eq:term1}
\end{equation}
and,

\begin{equation}
    g^q_2 = \frac{1}{1 - d_0 x_0} x_{j_1} \frac{1}{1 - d_1 x_0} x_{j_2} \ldots x_{j_q} \frac{1}{1 - d_q x_0} = g^{q-1}_2 \frac{x_{j_q}}{1 - d_q x_0}
\end{equation}
where $\{ p,q\} \in \mathbb{N} $ and $\{i_1,...,i_p,j_1,...,j_q\} \in \{0,1\}$. 

By assuming that the generating series is represented as a series of products, one can readily compute their shuffle products; results can be defined recursively as 
shuffle products of lower-order terms. This recursion ends when any of the following terms occur: $x \shuffle 1$, $1 \shuffle x$, or $1 \shuffle 1$. The process halts
as the shuffle product of a term with the identity element simply returns the term itself.

Looking in more detail: the term $\frac{1}{1-bx_0}$ can be rearranged into $1+\frac{bx_0}{1-bx_0}$; this operation can be applied to the highest-order fraction to give,

\begin{equation}
    g_1^p \shuffle g_2^q = \bigg[ g_1^{p-1}x_{i_p}\bigg(1+ \frac{b_p x_0}{1-b_p x_0} \bigg)\bigg] \shuffle \bigg[ g_2^{q-1}x_{j_q}\bigg(1+ \frac{d_q x_0}{1-d_q x_0} \bigg) \bigg]
\end{equation}

Recalling that the shuffle product is distributive over addition \cite{reutenauer1993freelie}; the above product can be expanded to give \cite{fliess1983algebraic},

\begin{multline}
    g_1^p \shuffle g_2^q  = \bigg[ g_1^{p-1} x_{i_p} \frac{1}{1 - b_p x_0} \shuffle g_2^{q-1} \bigg] x_{j_q} 
                          + \bigg[ g_1^{p-1} \shuffle g_2^{q-1} x_{j_q} \frac{1}{1 - d_q x_0} \bigg] x_{i_p} \\                      
                          + \bigg[ g_1^{p-1} x_{i_p} \frac{1}{1 - b_p x_0} \shuffle g_2^{q-1} x_{j_q} \frac{1}{1 - d_q x_0} \bigg] b_p x_0  
                          + \bigg[ g_1^{p-1} x_{i_P} \frac{1}{1 - b_p x_0} \shuffle g_2^{q-1} x_{j_q} \frac{1}{1 - d_q x_0} \bigg] d_q x_0
\end{multline}
(remembering that order of the terms is important).

Careful regrouping of terms gives a compact recursion for the shuffle product of two generating series of the form \cite{fliess1983algebraic},

\begin{equation}
    g_1^p \shuffle g_2^q = \frac{1}{1-(b_p+d_q)x_0}\bigg[ (g_1^p \shuffle g_2^{q-1})x_{j_q}+ (g_1^{p-1} \shuffle g_2^q)x_{i_p}\bigg]
\label{eq:shuff_rec}
\end{equation}

Noting that $\frac{1}{1 - b x_0} \shuffle \frac{1}{1 - d x_0} = \frac{1}{1 - (b + d) x_0}$.

Fortunately, the shuffle product is associative \cite{reutenauer1993freelie} i.e, $(g_{j_1}\shuffle g_{j_2})\shuffle g_{j_3}=g_{j_1}\shuffle (g_{j_2}\shuffle g_{j_3})$
and this means that equation (\ref{eq:shuff_rec}) can be extended straighforwardly (if tediously) to higher-order products, because the order in which the pairwise
products are taken does not matter.

One of the strengths of the generating series approach is that the formalism above is amenable to computer implementation; the calculations presented in this paper
are the result of a Python implementation of the necessary algebra. 

\subsection{The Asymmetric Duffing Oscillator}
\label{sec:gen_aduff}

When equation (\ref{eq:gen_poly}) is applied to the asymmetric Duffing equation in the canonical form in (\ref{eq:saduff}), the result is,

\begin{equation}
    g + x_0 g + x_0^2 g + \varepsilon_1 x_0^2 [g \shuffle g] + \varepsilon_2 x_0^2 [g \shuffle g \shuffle g] = x_0 x_1 
\label{eq:gen_saduff}
\end{equation}

By collecting like terms and factorising the quadratic expression in $x_0$, this equation can be rearranged into the form,

\begin{equation}
    g = \frac{x_0 x_1}{(1 - a_1 x_0)(1 - a_2 x_0)} - \varepsilon_1 \frac{x_0^2}{(1 - a_1 x_0)(1 - a_2 x_0)} g \shuffle g  
    - \varepsilon_2 \frac{x_0^2}{(1 - a_1 x_0)(1 - a_2 x_0)} g \shuffle g \shuffle g
\label{eq:fgen_saduff}
\end{equation}
where $(1 - a_1 x_0)(1 - a_2 x_0) = 1 + a x_0 + x_0^2$. 

The generating series solution can now be obtained for the oscillator; an iterative procedure can be followed where,

\begin{equation}
        g_{i+1} = - \frac{x_0}{1 - a_1 x_0}\frac{x_0}{1 - a_2 x_0} \bigg\{ \varepsilon_1 \sum_{i_1 + i_2 = i} [g_{i_1} \shuffle g_{i_2}] + 
        \varepsilon_2 \sum_{j_1 + j_2 + j_3 = i} [g_{j_1} \shuffle g_{j_2} \shuffle g_{j_3}] \bigg\}
\label{eq:gen_saduff_iter}
\end{equation}

The iteration begins with $g_1$ from the underlying linear system with $\varepsilon_1 = \varepsilon_2 = 0$,

\begin{equation}
     g_1 = \frac{1}{1 - a_1 x_0} x_0 \frac{1}{1 - a_2 x_0} x_1
\end{equation}

The generating series $g$, is then the sum,

\begin{equation}
    g = g_0 + g_1 + g_2 + \ldots + g_n + \ldots
\label{eq:gen_ser_saduff}
\end{equation}
allowing for the possibility for a constant offset $g_0$ in the response.

For the calculation here, only the first two iterations will be displayed. The additional nonlinear term here causes the iterations to have exponentially more terms 
compared to the single nonlinearities considered in other works \cite{fliess1983algebraic, fliess1982application}. 

For a more compact notation, the terms in the generating series, as shown in equation (\ref{eq:term1}) for example, can be expressed in the form of a $(2\times p)$ 
array; in the calculation here for example, one has,

\begin{table}[H]
\centering
$g_1 = 1 \begin{bmatrix} x_0  & x_1 \\ -a_1 & -a_2   \\ \end{bmatrix}$
\end{table}
\noindent where the notation shows a word in the numerator in the first row and the corresponding coefficient in the denominator in the second row. Each column 
represents a term in the rational fraction, and the overall coefficient/multiplier is found outside of the array.

The first iteration of the algorithm gives,

\begin{equation}
        g_{1} = - \frac{x_0}{1 - a_1 x_0} \frac{x_0}{1 - a_2 x_0} \big\{ \varepsilon_1 [g_{0} \shuffle g_{0}]
        + \varepsilon_2 [g_{0}\shuffle g_{0}\shuffle g_{0}] \big\}
\label{eq:iter_g1}
\end{equation}
(noting that order is unimportant in expressions considering a single letter).

Now expanding the shuffle products using equation (\ref{eq:shuff_rec}) yields,

\begin{table}[H]
\centering
     $g_1 = - 2 \varepsilon_1 \begin{bmatrix} x_0  & x_0 & x_0 & x_1  & x_0  & x_1 \\ 
                                              -a_1 & -a_2 & -2 a_1 & -a_1 - a_2 & -a_1 & -a_2 \\ \end{bmatrix}$
     $- 4 \varepsilon_1 \begin{bmatrix} x_0 & x_0 & x_0 & x_0 & x_1 & x_1 \\ 
                                      - a_1 & -a_2 & -2a_1 & -a_1 - a_2 & -2 a_2 & -a_2 \\ \end{bmatrix}$
     $- 6 \varepsilon_2 \begin{bmatrix} x_0 & x_0 & x_0 & x_1 & x_0 & x_1 & x_0  & x_1 \\
         - a_1 & -a_2 & -3 a_1 & -2 a_1 - a_2 & -2 a_1 & -a_1 - a_2 & - a_1 & - a_2 \\ \end{bmatrix}$
     $- 12 \varepsilon_2 \begin{bmatrix} x_0 & x_0 & x_0 & x_1    & x_0        & x_0    & x_1  & x_1 \\
                                       - a_1 & - a_2 & -3 a_1 & -2 a_1 - a_2 & -2 a_1 & -a_1 - a_2 & -2 a_2 & -a_2 \\ \end{bmatrix}$
     $- 24 \varepsilon_2 \begin{bmatrix} x_0 & x_0 & x_0 & x_0 & x_1 & x_0 & x_1 & x_1 \\
                                       - a_1 & -a_2 & -3 a_1 & -2 a_1 - a_2 & - a_1 - 2 a_2 & - a_1 - a_2 & - 2 a_2 & - a_2 \\ \end{bmatrix}$
     $- 12 \varepsilon_2 \begin{bmatrix} x_0 & x_0 & x_0 & x_0 & x_1 & x_1 & x_0 & x_1 \\
                                       - a_1 & - a_2 & - 3 a_1 & - 2 a_1 - a_2 & - a_1 - 2 a_2 & -a_1 - a_2 & - a_1 & - a_2 \\ \end{bmatrix}$
     $- 36 \varepsilon_2 \begin{bmatrix} x_0 & x_0 & x_0 & x_0 & x_0 & x_1 & x_1 & x_1 \\
                                       - a_1 & -a_2 & - 3 a_1 & - 2 a_1 - a_2 & - a_1 - 2 a_2 & - 3 a_2 & - 2 a_2 & -a_2 \\ \end{bmatrix}$
\end{table}
and for $g_2$,

\begin{equation}
        g_{2} = - \frac{x_0}{1 - a_1 x_0}\frac{x_0}{1 - a_2 x_0} \big\{ \varepsilon_1 [g_{1} \shuffle g_{0} + g_{0} \shuffle g_{1}] + 
        \varepsilon_2 [g_{1} \shuffle g_{0} \shuffle g_{0} + g_{0} \shuffle g_{1} \shuffle g_{0} + g_{0} \shuffle g_{0} \shuffle g_{1}] \big\}
\label{eq:iter_g2}
\end{equation}

\begin{table}[H]
\centering
      $ g_{2} = 6 \varepsilon_1^2 \begin{bmatrix} x_0 & x_0 & x_0 & x_0 & x_0 & x_1 & x_0 & x_1 & x_0  & x_1 \\
                                                - a_1 & - a_2 & - 2 a_1 & - a_1 - a_2 & -3 a_1 & -2 a_1 - a_2 & -2 a_1 & -a_1 - a_2 & -a_1 & -a_2 \\ \end{bmatrix}$
      $+ 8 \varepsilon_1^2 \begin{bmatrix} x_0 & x_0 & x_0 & x_0 & x_0 & x_1 & x_0 & x_0 & x_1 & x_1 \\
                                         - a_1 & - a_2 & -2 a_1 & - a_1 - a_2 & -3 a_1 & -2 a_1 - a_2 & -2 a_1 & -a_1 - a_2 & -2 a_2 & -a_2 \\ \end{bmatrix}$
      $+ 8 \varepsilon_1^2 \begin{bmatrix} x_0 & x_0 & x_0 & x_0 & x_0 & x_0 & x_1 & x_0 & x_1 & x_1 \\
                                         - a_1 & - a_2 & -2 a_1 & - a_1 - a_2 & -3 a_1 & -2 a_1 - a_2 & - a_1 - 2 a_2 & - a_1 - a_2 & -2 a_2 & - a_2 \\ \end{bmatrix}$
      $+ 8 \varepsilon_1^2 \begin{bmatrix} x_0 & x_0 & x_0 & x_0 & x_0 & x_0 & x_1 & x_0 & x_1 & x_1 \\
                                         - a_1 & - a_2 & -2 a_1 & - a_1 - a_2 & -2 a_2 & -2 a_1 - a_2 & - a_1 - 2 a_2 & - a_1 - a_2 & -2 a_2 & - a_2 \\ \end{bmatrix}$
      $+ 8 \varepsilon_1^2 \begin{bmatrix} x_0 & x_0 & x_0 & x_0 & x_0 & x_0 & x_1 & x_1 & x_0 & x_1 \\
                                         - a_1 & - a_2 & -2 a_1 & - a_1 - a_2 & -3 a_1 & -2 a_1 - a_2 & - a_1 - 2 a_2 & - a_1 - a_2 & - a_1 & -a_2 \\ \end{bmatrix}$
      $+ 8 \varepsilon_1^2 \begin{bmatrix} x_0 & x_0 & x_0 & x_0 & x_0 & x_0 & x_1 & x_1 & x_0 & x_1 \\
                                         - a_1 & -a_2 & -2 a_1 & - a_1 - a_2 & -2 a_2 & -2 a_1 - a_2 & - a_1 - 2 a_2 & - a_1 - a_2 & - a_1 & - a_2 \\ \end{bmatrix}$
      $+ 4 \varepsilon_1^2 \begin{bmatrix} x_0 & x_0 & x_0 & x_0 & x_0 & x_1 & x_0 & x_1 & x_0 & x_1 \\
                                         - a_1 & - a_2 & -2 a_1 & - a_1 - a_2 & -2 a_2 & -2 a_1 - a_2 & -2 a_1 & - a_1 - a_2 & - a_1 & - a_2 \\ \end{bmatrix}$
      $+ 4 \varepsilon_1^2 \begin{bmatrix} x_0 & x_0 & x_0 & x_0 & x_1 & x_0 & x_0 & x_1 & x_0 & x_1 \\
                                         - a_1 & - a_2 & -2 a_1 & - a_1 - a_2 & -2 a_2 & - a_2 & -2 a_1 & - a_1 - a_2 & - a_1 & - a_2 \\ \end{bmatrix}$
      $+ 2 \varepsilon_1^2 \begin{bmatrix} x_0 & x_0 & x_0 & x_1 & x_0 & x_0 & x_0 & x_1 & x_0 & x_1 \\
                                         - a_1 & - a_2 & -2 a_1 & - a_1 - a_2 & - a_1 & - a_2 & -2 a_1 & - a_1 - a_2 & - a_1 & - a_2    \\ \end{bmatrix}$
      $+ 4 \varepsilon_1^2 \begin{bmatrix} x_0 & x_0 & x_0 & x_0 & x_0 & x_0 & x_1 & x_1 & x_0 & x_1 \\
                                         - a_1 & - a_2 & -2 a_1 & - a_1 - a_2 & -3 a_1 & -2 a_1 - a_2 & - a_1 - 2 a_2 & - a_1 - a_2 & - a_1 & - a_2 \\ \end{bmatrix}$
\end{table}

Only the first 10 terms have been shown for $g_2$, as there are 360 terms in the full expansion.

\section{Determining System Response} 
\label{sec:response}

The analysis up to now has allowed the input-output relationship for the system to be expressed in terms of the generating series algebra. In order to compute an actual
response, one needs to substitute for the relevant excitation $x(t)$, as encoded in the letter $x_1$ in the free monoid, and then transform back to the time domain. The
analysis has provided terms of a specific form; words in the letters $x_0$, and $x_1$ and rational fractions of them. If the transformation back is made with general $x_1$,
the result will be a Volterra expansion, and one will be able to read off the Volterra kernels. Each term in the series -- of the form given by equation (\ref{eq:gen_fact2})
-- corresponds to a specific iterated integral. Each appearance of $x_0$ represents an integration, so a term with $g$ occurrences of $x_0$ represents a $q$-fold iterated
integral. Fliess and co-workers computed the general inverse transform of a $q$-fold product in the generating series, it takes the form \cite{fliess1983algebraic},

\begin{equation}
    \int_0^t \int_0^{\tau_q} \ldots \int_0^{\tau_2} f_{a_1}^{\alpha_1}(t - \tau_q) \ldots f_{a_{q-1}}^{\alpha_{q-1}}(\tau_2 - \tau_1) f_{a_{q}}^{\alpha_{q}}(\tau_1) x(\tau_q) \ldots x(\tau_1) d\tau_q \ldots \tau_1
\label{eq:it_volt}
\end{equation}
where
\begin{equation}
    f_a^\alpha = \bigg[ \sum_{j=0}^{\alpha - 1} {\alpha - 1 \choose j} \frac{a^j t^j}{j!} \bigg] e^{at}
\end{equation}

This relationship shows how the generating series and Volterra series are so strongly linked. The correspondence between the terms in equation (\ref{eq:gen_ser_saduff}) and 
the Volterra terms can be shown to be,

\begin{align}
   g_0  &\Leftrightarrow  y_0=h_0\\
   g_1  &\Leftrightarrow y_1=\int_{-\infty}^{+\infty}h_1(\tau_1)x(\tau_1)d\tau_1 \\
   g_2  &\Leftrightarrow y_2=\int_{-\infty}^{+\infty}\int_{-\infty}^{+\infty} h_2(\tau_2,\tau_1)x(\tau_2)x(\tau_1)d\tau_2d\tau_1
\end{align}

By computing the inverse Laplace-Borel transform of the generating series derived for the oscillators in Section \ref{sec:gen_aduff}, the Volterra kernels are determined.
For specific excitations $x(t)$, the system response can be computed. The relevant inverse transforms are tabulated below:

\begin{table}[H]
    \centering
\begin{tabular}{|c|c|}
\hline
$x(t)$           & $g[x(t)]$               \\ \hline
Unit Step        & $1$                        \\
$\frac{t^n}{n!}$ & $x_0^n$                  \\
$\bigg( \sum_{i=0}^{n-1} \binom{i}{n-1}\frac{a^it^i}{i!} \bigg) e^{at}$              & $(1-ax_0)^{-n}$          \\
$\cos(\omega t)$  & $(1+\omega^2x_0^2)^{-1}$ \\ \hline
\end{tabular}
\caption{Laplace-Borel Transforms of Common Functions \cite{li2004formal,fliess1983algebraic}}
\label{tab:lap bor}
\end{table}

Rather than give the asymmetric Duffing system response for a deterministic excitation, a little more work will allow characterisation of the response to a random
excitation.

\section{Response to Gaussian White Noise}
\label{sec:noise}

Clearly, the machinery provided up to now can not provide a time-series response to a truly random excitation; however, it can be adapted to give output statistics,
of the response and this is simplest when the excitation is a {\em Gaussian white noise} process. Such a process is specified by a demand that it have zero mean, and an 
auto-correlation function of the form \cite{arnold1974stochastic},

\begin{equation}
    E[x(t)x(\tau)] = \langle x(t) x(\tau) \rangle = \sigma^2 \delta (t - \tau)
\label{eq:auto}
\end{equation}
where $\sigma^2$ denotes the `noise-power' (one must take care in interpreting this as a variance) and angle brackets denote expectations. 

The most basic statistic one can estimate is the mean of the response, or its expectation $\langle y(t) \rangle$. In \cite{bedrosian1971output}, the authors developed an appropriate form of the Volterra series for random excitation, based on stochastic calculus \cite{ito1944}. This formulation was adapted by Fliess \cite{fliess1982application},
in order to compute statistics from the generating series. The expectations are interpreted as ensemble averages so that one can take the expectations in the transform
domain and then map back. In this way $\langle g \rangle$ corresponds to $\langle y(t) \rangle$, and the higher-order statistics $\langle y(t)^n \rangle$ are obtained
by mapping back the shuffle products $\langle g^{\shuffle n} \rangle$.  

As usual now, it is sufficient to consider only the calculation for terms of the form shown in equation (\ref{eq:gen_fact2}); the basic rules are \cite{lamnabhi1983stochasticcomp},

\begin{equation}
     \big\langle \frac{1}{1 - b_0 x_0} x_{i_1} \frac{1}{1 - b_1 x_0} \ldots  x_{i_n} \frac{1}{1 - b_n x_0} \big\rangle =
\begin{dcases}
     \frac{x_{0}}{1 - b_0 x_0} \big\langle \frac{1}{1 - b_1 x_0} x_{i_2} \ldots x_{i_n} \frac{1}{1 - b_n x_0} \big\rangle, & \text{if} \ i_1=0 \\
     \frac{\sigma^2}{2} \frac{x_{0}}{1 - b_0 x_0} \big\langle \frac{1}{1 - b_2 x_0} x_{i_3} \ldots x_{i_n} \frac{1}{1 - b_n x_0} \big\rangle, & \text{if} \ i_1=i_2=1 \\
     0, & \text{Otherwise}
\label{eq:cond}
\end{dcases}
\end{equation}

Note how restrictive this recursion is, many terms will automatically be zero; this is related to the fact that expectations of products of an odd number of Gaussian
random variables will average to zero. Once the generating series has been decomposed into the standard terms (and in this case will only contain the letter $x_0$), 
the usual rules allow inversion using a table of Laplace-Borel transforms.

The autocorrelation of the response is a little more complicated, this has the form,

\begin{equation}
     S_{yy}= \langle y(t_1)y(t_2)] \rangle
\label{eq:resp_auto}
\end{equation}
and the product of $y$s will produce a shuffle product in the domain of the generating series.

For the asymmetric Duffing oscillator under investigation here, performing the ensemble average for all the terms in the generating series, the following result is obtained,

\begin{multline}
\langle g \rangle = -4\varepsilon_1 \left(\frac{\sigma^2}{2}\right)
\begin{bmatrix}
  x_0  & x_0    & x_0        & x_0    & x_0  \\
 -a_1     & -a_2 & -2a_1 & -a_1 - a_2 & -2a_2    \\
\end{bmatrix}\\
+48 \varepsilon_1 \varepsilon_2 \left(\frac{\sigma^2}{2}\right)^2
\begin{bmatrix}
  x_0 & x_0 & x_0 & x_0 & x_0 & x_0 & x_0 & x_0 & x_0 & x_0  \\
 -a_1 & -a_2 & -2a_1 & -a_1 - a_2 & -4a_1 & -3a_1-a_2 & -2a_1-2a_2 & 2a_1& -a_1-a_2 & -2a_2    \\
\end{bmatrix}\\
+48 \varepsilon_1 \varepsilon_2 \left(\frac{\sigma^2}{2}\right)^2
\begin{bmatrix}
  x_0 & x_0 & x_0 & x_0 & x_0 & x_0 & x_0 & x_0 & x_0 & x_0  \\
 -a_1 & -a_2 & -2a_1 & -a_1 - a_2 & -2a_1 & -3a_1-a_2 & -2a_1-2a_2 & 2a_1& -a_1-a_2 & -2a_2    \\
\end{bmatrix}\\
+ \ldots
\label{eq:gwn_resp}
\end{multline}

As discussed above, because of the third condition in the recursion (\ref{eq:cond}), many terms are zero and do not contribute. For example, the second significant term 
in the expansion above is actually the 39$^{th}$ term from the total of 360 in $g_2$. 

To move the computation forward, it is necessary to carry out the partial fraction calculations implicit in the terms in equation (\ref{eq:gwn_resp}). At this point,
it is useful to specify numerical values for $a_1$ and $a_2$, as the partial fractions calculations can be cumbersome when carried out algebraically. In the case of the
asymmetric Duffing equation, the system parameters are $m, c, k_1, k_2, k_3$. The relevant parameters in the generating series can then be calculated via the scaled 
version of the Duffing equation in equation(\ref{eq:saduff}); starting with values here of $m =1kg$, $c = 15Nsm^{-1}$, $k_1 = 25Nm^{-1}$, $k_2 = 625Nm^{-1}$ and 
$k_3 = 7500Nm^{-1}$, one obtains $a_1 = -\frac{1}{2}(3+\sqrt{5})$, $a_2 = -\frac{1}{2}(3 + \sqrt{5})$, $\varepsilon_1 = 1$ and $\varepsilon_2 = 0.5$. 

By decomposing the terms in equation (\ref{eq:gwn_resp}) into partial fractions and applying the inverse Laplace-Borel transforms as given in Table 1, the mean 
response of the system $E[y(t)]$ can be computed; the result is,

\begin{multline}
    \langle y(t) \rangle =
    -4\varepsilon_1 \left( \frac{\sigma^2}{2}\right)\bigg[ 0.08332 -\frac{0.0005835}{0.1910}e^{-\frac{t}{0.1910}} +\frac{0.02288}{0.3333}e^{-\frac{t}{0.3333}} - \frac{0.6315}{2.618}e^{-\frac{t}{2.618}} \\ - \frac{0.03514}{0.3820}e^{-\frac{t}{0.3820}} + \frac{0.2409}{1.309}e^{-\frac{t}{1.309}} \bigg]
    +48\varepsilon_1 \varepsilon_2 \left( \frac{\sigma^2}{2}\right)\bigg[ 0.0001508+ \frac{0.0005778}{0.01667}e^{-\frac{t}{0.01667}} \\ - \frac{0.001053}{1.309^2}(1-\frac{t}{1.309})e^{-\frac{t}{1.309}}) - \frac{0.0002746}{0.1910}e^{-\frac{t}{0.1910}}  +  \frac{0.000006524}{0.09551}e^{-\frac{t}{0.09551}} + \frac{0.0005542}{0.3333}e^{-\frac{t}{0.3333}} \\ - \frac{0.003057}{2.618}e^{-\frac{t}{2.618}} - \frac{0.00002157}{0.1910^2}(1-\frac{t}{0.1910})e^{-\frac{t}{0.1910}}  - \frac{0.00005400}{0.1214}e^{-\frac{t}{0.1214}} \\ -  \frac{0.0007346}{0.3820}e^{-\frac{t}{0.3820}} - \frac{0.0007298}{0.3333^2}(1-\frac{t}{0.3333})e^{-\frac{t}{0.3333}} + \frac{0.001947}{1.309}e^{-\frac{t}{1.309}} \bigg] + \ldots
\end{multline}

\begin{figure}[htbp!]
\centering
\includegraphics[width=0.6\textwidth]{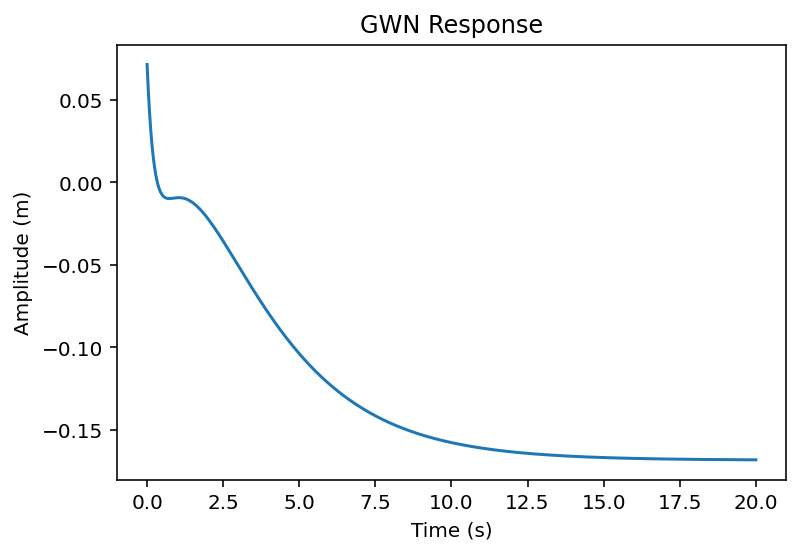}
\caption{Expectation of the response of an asymmetric Duffing oscillator to a Gaussian white-noise excitation.}
\label{fig:GWN_response}
\end{figure}

Choosing the somewhat arbitrary value $\sigma=1$, the result in Figure \ref{fig:GWN_response} is obtained. One observes a transient which occurs from `switching on' the
excitation at $t = 0$. In fact, because the asymmetric Duffing oscillator has stationary response if the input is stationary, the expectation will tend to a constant 
value as $t \longrightarrow \infty$; because the restoring force is asymmetric, that constant value will be non-zero.

\section{Conclusions}
\label{sec:conclusions}

Long conclusions are not warranted here as the aim of this paper was simply to revisit the generating series of Fliess and co-workers, as an elegant means of nonlinear
system analysis. In order to introduce a novel element, the analysis has been extended beyond previous work in order to deal with the case of two nonlinear terms 
in the equation of motion. This new analysis is also extended to the diagrammatic representation, where the presence of two nonlinearities produces two types of 
vertices in the `Feynman' rules for the diagrams.

Further work on the series is considering how it can be used in an automated manner in order determine Higher-order Frequency Response Functions for nonlinear 
structural dynamic systems.

\section*{Acknowledgements}

The authors would like to thank the UK EPSRC for funding through the Established Career Fellowship EP/R003645/1 and the Programme Grant EP/R006768/1.

\bibliographystyle{unsrt}
\bibliography{imac_21_TG_1}

\begin{thebibliography}{10}

\bibitem{worden2019nonlinearity}
K.\ Worden and G.R.\ Tomlinson.
\newblock {\em Nonlinearity in Structural Dynamics: Detection, Identification
  and Modelling}.
\newblock Institute of Physics Publishing, 2001.

\bibitem{volterra1930theory}
V.\ Volterra.
\newblock {\em Theory of Functionals and of Integral and Integro-Differential
  Equations}.
\newblock Blackie \& Son Limited, 1930.

\bibitem{barrett1963use}
J.F.\ Barrett.
\newblock The use of functionals in the analysis of non-linear physical
  systems.
\newblock {\em International Journal of Electronics}, 15:567--615, 1963.

\bibitem{fliess1981fonctionnelles}
M.\ Fliess.
\newblock Fonctionnelles causales non lin{\'e}aires et ind{\'e}termin{\'e}es
  non commutatives.
\newblock {\em Bulletin de la Soci{\'e}t{\'e} Math{\'e}matique de France},
  109:3--40, 1981.

\bibitem{fliess1982application}
M.\ Fliess and F.\ Lamnabhi-Lagarrigue.
\newblock Application of a new functional expansion to the cubic anharmonic
  oscillator.
\newblock {\em Journal of Mathematical Physics}, 23:495--502, 1982.

\bibitem{lamnabhi1982new}
M.\ Lamnabhi.
\newblock A new symbolic calculus for the response of nonlinear systems.
\newblock {\em Systems \& Control Letters}, 2:154--162, 1982.

\bibitem{fliess1983algebraic}
M.\ Fliess, M.\ Lamnabhi, and F.\ Lamnabhi-Lagarrigue.
\newblock An algebraic approach to nonlinear functional expansions.
\newblock {\em IEEE Transactions on Circuits and Systems}, 30:554--570, 1983.

\bibitem{lamnabhi1986functional}
M.\ Lamnabhi.
\newblock Functional analysis of nonlinear circuits: a generating power series
  approach.
\newblock {\em IEE Proceedings H (Microwaves, Antennas and Propagation)},
  133:375--384, 1986.

\bibitem{lamnabhi1980application}
F.\ Lamnabhi-Lagarrigue.
\newblock {\em Application des variables non commutatives {\`a} des calculs
  formels en statistique non lin{\'e}aire}.
\newblock PhD thesis, Universit{\'e} Paris-Sud, 1980.

\bibitem{lamnabhi1982nonlinearcomp}
F.\ Lamnabhi-Lagarrigue and M.\ Lamnabhi.
\newblock Algebraic computation of the solution of some nonlinear differential
  equations.
\newblock In {\em Proceedings of the European Computer Algebra Conference},
  pages 204--211. Springer, 1982.

\bibitem{li2004formal}
Y.\ Li and W.S.\ Gray.
\newblock The formal {L}aplace-{B}orel transform, {F}liess operators and the
  composition product.
\newblock In {\em Proceedings of the 36$^{th}$ Southeastern Symposium on System
  Theory}, pages 333--337. IEEE, 2004.

\bibitem{duffing1918erzwungene}
G.\ Duffing.
\newblock {\em Erzwungene Schwingungen bei Ver{\"a}nderlicher Eigenfrequenz und
  ihre Technische Bedeutung}.
\newblock F.\ Vieweg and Sohn, 1918.

\bibitem{reutenauer1993freelie}
C.\ Reutenauer.
\newblock {\em Free Lie Algebras}.
\newblock Elsevier, 1993.

\bibitem{ree1958lie}
R.\ Ree.
\newblock Lie elements and an algebra associated with shuffles.
\newblock {\em Annals of Mathematics}, pages 210--220, 1958.

\bibitem{bedrosian1971output}
E.\ Bedrosian and S.O.\ Rice.
\newblock The output properties of {V}olterra systems (nonlinear systems with
  memory) driven by harmonic and {G}aussian inputs.
\newblock {\em Proceedings of the IEEE}, 59:1688--1707, 1971.

\bibitem{lubbock1969multidimensional}
J.K.\ Lubbock and V.S.\ Bansal.
\newblock Multidimensional {L}aplace transforms for solution of nonlinear
  equations.
\newblock {\em Proceedings of the Institution of Electrical Engineers},
  116:2075--2082, 1969.

\bibitem{morton1970}
J.B.\ Morton and S.\ Corrsin.
\newblock Consolidated expansions for estimating the response of a randomly
  driven nonlinear oscillator.
\newblock {\em Journal of Statistical Physics}, 2:153--194, 1970.

\bibitem{feynman1948space}
R.P.\ Feynman.
\newblock Space-time approach to non-relativistic quantum mechanics.
\newblock {\em Reviews of Modern Physics}, 20:367--387, 1948.

\bibitem{gantmakher2000theory}
F.R.\ Gantmakher.
\newblock {\em The Theory of Matrices}.
\newblock American Mathematical Society, 2000.

\bibitem{arnold1974stochastic}
L.\ Arnold.
\newblock {\em Stochastic Differential Equations}.
\newblock Wiley-Blackwell, 1974.

\bibitem{ito1944}
K.\ Ito.
\newblock Stochastic integral.
\newblock {\em Proceedings of the Imperial Academy}, 20:519--524, 1944.

\bibitem{lamnabhi1983stochasticcomp}
F.\ Lamnabhi-Lagarrigue and M.\ Lamnabhi.
\newblock Algebraic computation of the statistics of the solution of some
  nonlinear stochastic differential equations.
\newblock In {\em Proceedings of the European Conference on Computer Algebra},
  pages 55--67. Springer, 1983.

\end{thebibliography}

\end{document}